\documentclass[a4paper,10pt,twoside]{article}
\usepackage{RAstyle}
\usepackage[german,czech,english]{babel}

\usepackage[latin2]{inputenc}    
\usepackage{epsfig}               
\usepackage{epstopdf}
\usepackage[small,bf]{caption}
\DeclareCaptionLabelFormat{fig}{\hspace{2.5mm}Fig. #2}
\DeclareCaptionLabelFormat{tab}{\hspace{2.5mm}Tab. #2}
\DeclareCaptionLabelSeparator{dot}{. }
\newcommand{\fcaption}[1]{
\captionsetup{labelformat=fig}
\caption{#1}
}

\usepackage{amssymb,amsmath}
\newcommand{\ZZ}{\mathbb{Z}}
\newcommand{\RR}{\mathbb{R}}
\newcommand{\CC}{\mathbb{C}}

\newcommand{\re}{\mathrm{Re\,}}

\begin{document}


\twocolumn[
\begin{center}{
	{\huge \bf Trigonometric Splines for Oscillator Simulation} \\
	\vspace{6mm}
	{\it \large Kai BITTNER\footnotemark[1], Hans-Georg BRACHTENDORF\footnotemark[1]\\}
	\vspace{4mm}
	{\footnotemark[1] University of Applied Sciences Upper Austria, Softwarepark 11, 4232 Hagenberg, Austria}\\
	\vspace{4mm}
	{Kai.Bittner@fh-hagenberg.at,\quad Hans-Georg.Brachtendorf@fh-hagenberg.at}\\
	\ \\
	\ \\}
\end{center}
]


\selectlanguage{english} 

\abstract{We investigate the effects of numerical damping for oscillator simulation with spline methods.
Numerical damping results in an artificial loss of energy and leads therefore to unreliable results
in the simulation of autonomous systems, as e.g.\   oscillators. We show that the negative effects of numerical damping
can be eliminated by the use of trigonometric splines. This will be in particular important for spline based adaptive methods. 
}
\vspace{11mm}



\section{Introduction}\label{sec:name1}

The simulation of oscillators suffers often from a loss of energy due to numerical damping. On the other hand 
numerical damping is often required to cancel out oscillations which occur due to numerical noise.
In the recent years, the expansion of the signal waveforms by wavelets or splines was investigated
\cite{SN03,SGN07,DCB05,Dau05,Bra2009,BiDau10b, BiDau10a}. This has been motivated by the fact
that trigonometric basis  (originally used in many RF-simulations) are not suited for the representation
of pulse-shaped signals due to a slow convergence
and the Gibbs' phenomenon.

\begin{figure}[h]
	\center
	\includegraphics[width=\columnwidth]{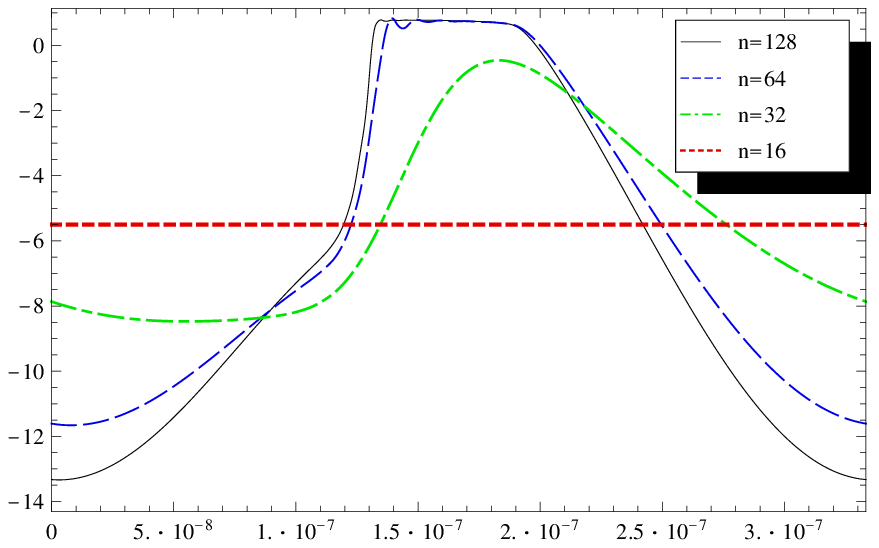}
	\fcaption{Spline solutions for a 3MHz Colpitz-Quartz oscillator.}\label{f:fig1}
\end{figure}

Here, we want to investigate the effects of numerical damping on a oscillator simulation by a spline collocation method
on uniform grids. Numerical experiments suggest that the effect described and investigated here do also occur
on nonuniform grids, which may appear in an adaptive method. However, we restrict to uniform grids, since this 
restriction permits the theoretical understanding of the observed effects. 

\begin{figure}[h]
	\center
	\includegraphics[width=\columnwidth]{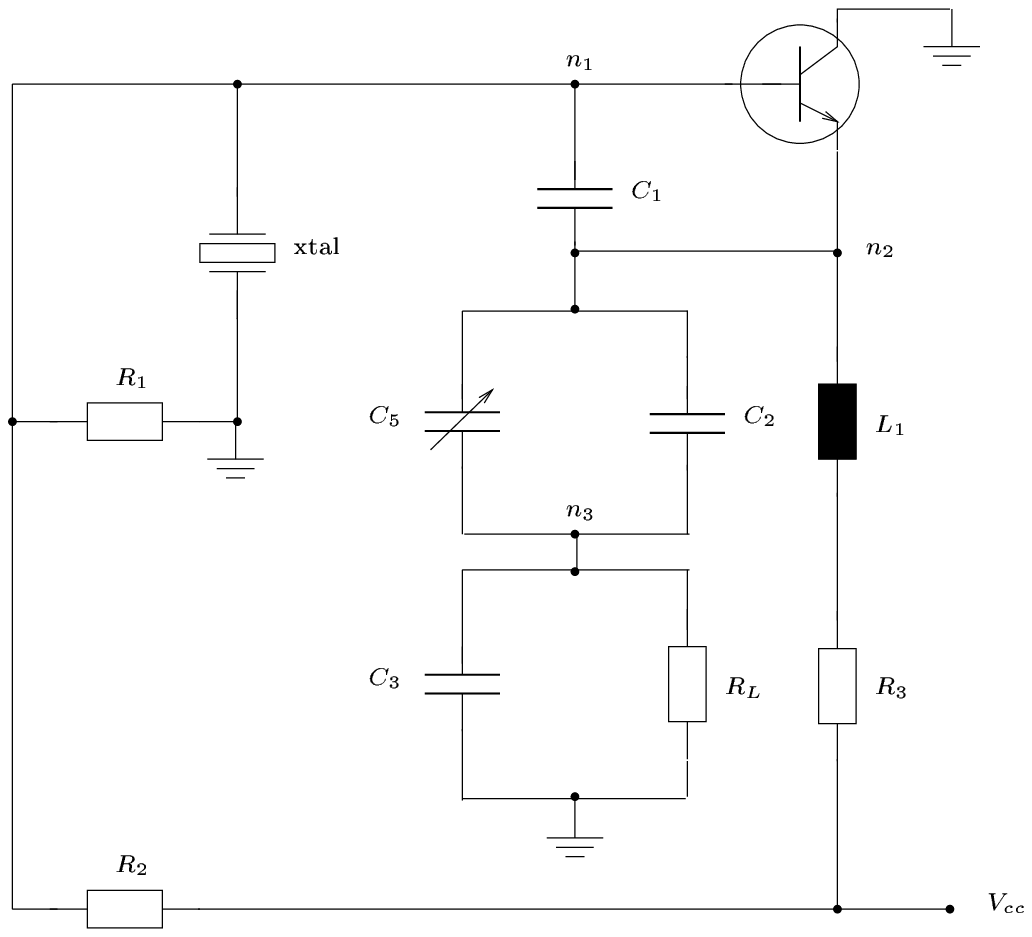}
	\fcaption{Schematic of the 3MHz Colpitz-Quartz oscillator.}\label{f:fig1a}
\end{figure}

In Fig.~\ref{f:fig1} we see numerical solutions for the periodic steady state of a  3MHz Colpitz-Quartz oscillator 
(see Fig.~\ref{f:fig1a})
computed by a spline collocation method on various grids. Although one expects an increase of the approximation error
on coarser grids, we can observe here an additional effect, which is a smaller amplitude for the coarser grid.
This artificial loss of energy is due to an effect called numerical damping, which is caused not by the physical system
or its mathematical model, but  by properties of the numerical method. 

In Sect.~\ref{sec:damping} we will investigate, the numerical damping in a spline collocation method.
We will overcome the negative effects of this numerical damping by the application of trigonometric splines
in Sect.~\ref{sec:trig}.

\pagebreak

\section{Numerical Damping for Splines}\label{sec:damping}

We consider the cardinal B-spline $N_m$ of order $m$ defined by the recursion
\begin{eqnarray*}
N_1(t) &:=& \chi_{(0,1]}(t) :=
\left\{\begin{array}{ll}1,&\mbox{for } t\in(0,1],\\
0&\mbox{otherwise,}\end{array}\right.\\[1ex]
N_m(t) &:=& \frac{t\,N_{m-1}(t)+(m-t)\,N_{m-1}(t-1)}{m-1}.
\end{eqnarray*}
For a detailed introduction to splines we refer to \cite{Sch81}.
It is well known that the family $\{N_m(\frac{t}{h}-k):~k\in\ZZ\}$
constitutes a stable basis for the spline space
$$
S_{m,h}=\Big\{f\in C^{m-2}(\RR):~f\big|_{\big(h k,h(k+1)\big)}\in\Pi_{m-1}\}
$$
of piecewise polynomials of degree less than $m$, which are $m-2$ times differentiable. 
Here, $\Pi_n$ denotes the space of polynomials of degree up to $n$

We consider now how the differential operator is approximated by a spline collocation method.
For the simplicity of the presentation and without loss of generality we restrict ourselves to 1-periodic functions (i.e.~with unity period).
The statements can be generalized to other period lengths by scaling.
First the function $x(t)$ is interpolated at the collocation points  $t_k=\frac{k+\frac{m}{2}+\sigma}{n}$
by a spline $s\in S_{m,\frac{1}{n}}$,
i.e., we have to determine coefficients  $c_\ell$ such that
\begin{equation}\label{interpol}
y_k:=x(t_k)=s(t_k)=\sum_\ell c_\ell\, N_m(k-\ell+\tfrac{m}{2}+\sigma).
\end{equation}
Here, the parameter $\sigma$ describes the deviation of the collocation points from
the center of the B-splines in relation to the mesh size.
The choice of $\sigma$ influences the stability of the spline interpolation, but also
properties of a collocation scheme for differential equations, as we will see in the sequel.

Applying the discrete Fourier transform
$$
\hat{x}_k = \sum_{\ell=0}^{n-1} x_\ell e^{\frac{2\pi i k\ell}{n}}
$$
to (\ref{interpol}) we obtain
$$ 
\hat{c}_k= \frac{\hat{y}_k}{\phi_m(\sigma,\frac{k}{n})},
$$
where 
$$
\phi_m(x,\xi)=\sum_{k\in\ZZ}N_m(x+\tfrac{m}{2}+k)\,e^{2\pi i k \xi}
$$
is Schoenberg's exponential Euler spline \cite{JRS91}.
Analogously, we obtain the values of the derivatives as
$$
y'_k:=s'(t_k)=\sum_\ell c_\ell N^\prime_m(k-\ell+\tfrac{m}{2}+\sigma),
$$
which yields 
\begin{equation}\label{spline_diff}
\widehat{y^\prime_k}=\psi_m(\sigma,\tfrac{k}{n})\;\hat{y}_k,
\end{equation}
with
$$
\psi_m(x,\xi)=\frac{\tfrac{\partial}{\partial x}\phi_m(x,\xi)}{\phi_m(x,\xi)}.
$$
Obviously $\psi_m(x,\xi)$ is 1-periodic in both arguments.

\begin{figure}[ht]
	\center
	\includegraphics[width=.9\columnwidth]{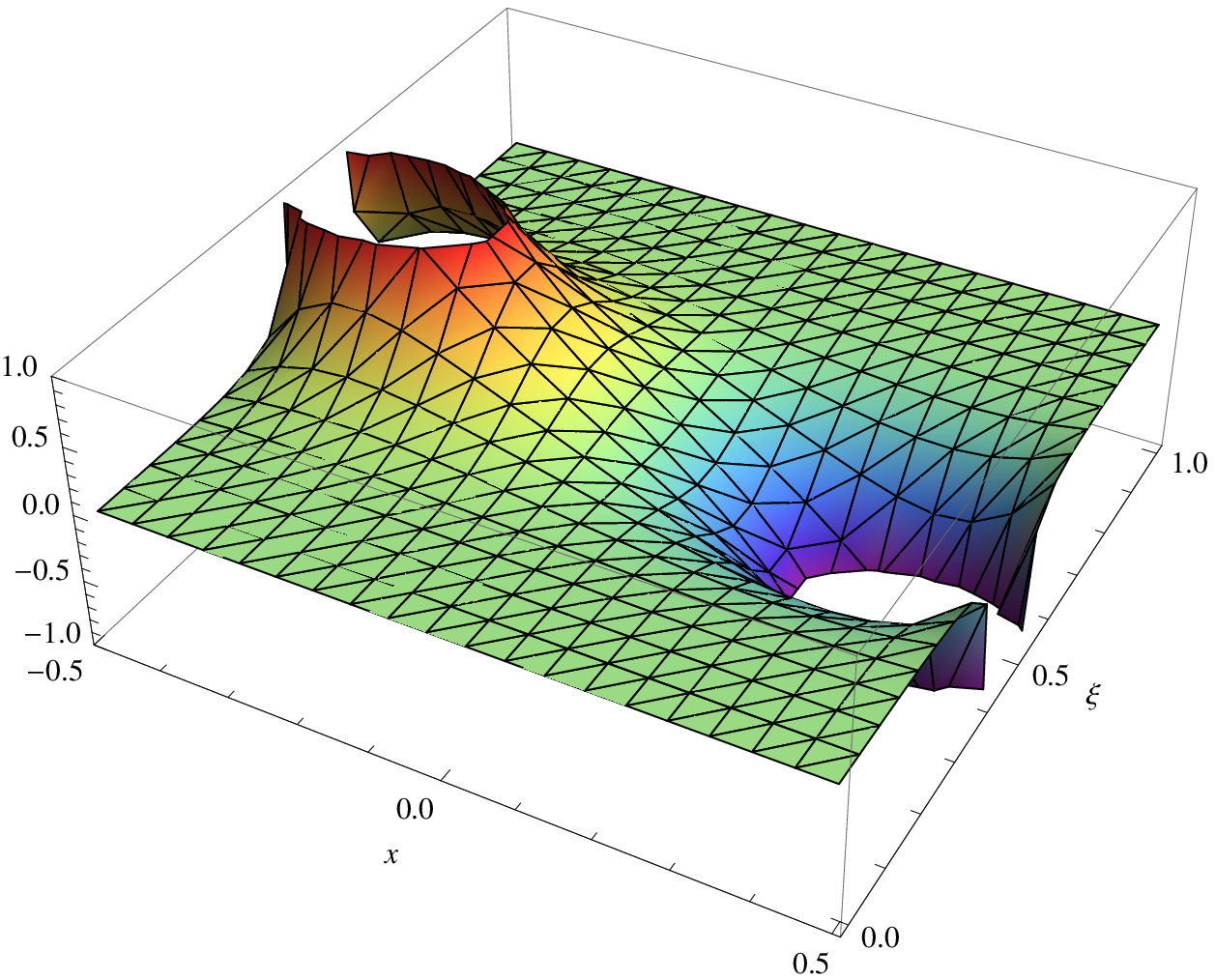}
	\fcaption{Plot of $\re \psi_3(x,\xi)$.}\label{f:fig2}
\end{figure}

Fig.~\ref{f:fig2} shows a plot of the real part of $\psi_3(x,\xi)$.
Due to a zero of $\phi_m$ at $(\frac{1}{2},\frac{1}{2})$ we have a singularity
of $\psi_m(x,\xi)$ at this point. This means that the numerical derivative becomes instable for 
$|\sigma|\approx\frac{1}{2}$ and we restrict ourselves to $|\sigma|<\frac{1}{2}$ sufficiently small. 
Furthermore, $\phi_m(x,0) = 1$, which implies  $\psi_m(x,0) = 0$.
Since $N_m(x+\frac{m}{2})$ is even we conclude that $\phi_m(0,\xi)$ is real, while $\tfrac{\partial}{\partial x}\phi_m(0,\xi)$
is purely imaginary, i.e., $\re \psi_m(0,\xi)= 0$. 

\begin{figure}[h]
	\center
	\includegraphics[width=\columnwidth]{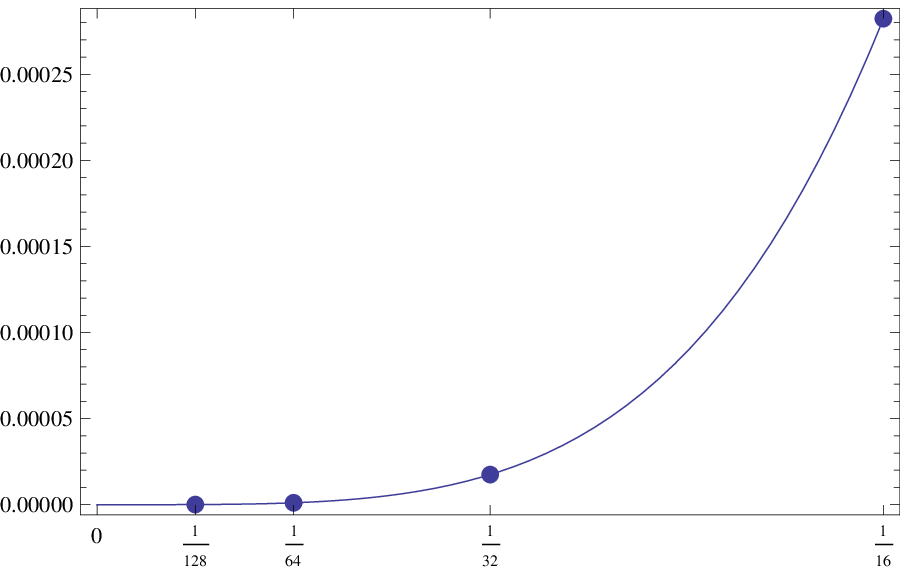}
	\fcaption{Plot of $\re \psi_3(-\frac{1}{4},\xi)$.}\label{f:fig3}
\end{figure}

Apparently $\re \psi_m(\sigma,\xi) > 0$ for $\sigma\in(-\frac{1}{2},0)$, which causes a numerical damping. 
For the
fundamental frequency the damping depends on the size of $\re \psi_m(\sigma,\frac{1}{n}) > 0$, where $n$ is the grid size.
In Fig.~\ref{f:fig3} we see the corresponding values for the simulations in  Fig.~\ref{f:fig1}, where $\sigma=-\frac{1}{4}$
was used. Obviously the grid size $n$ has to be chosen sufficiently large, in order to avoid a loss of energy by numerical damping.

\begin{figure}[ht]
	\center
	\includegraphics[width=\columnwidth]{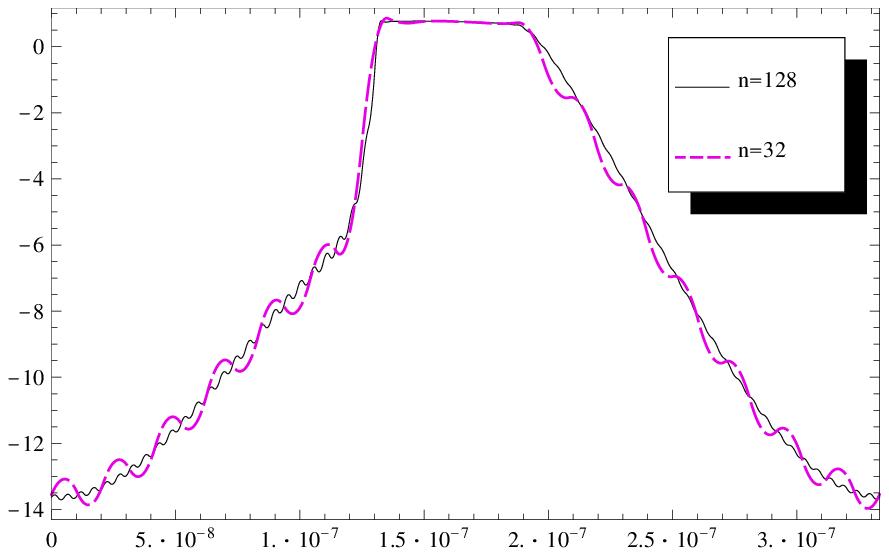}
	\fcaption{Spline solutions for the 3MHz  Colpitz-Quartz oscillator, with $\sigma=-.01$.}\label{f:fig4}
\end{figure}

The numerical damping is reduced for $\sigma\to 0$, $\sigma>0$.\linebreak This can be seen in  Fig.~\ref{f:fig4}, where almost no damping
of the fundamental frequency can be observed. However, there occurs also no damping of the high frequencies, which is
necessary in a simulation to eliminate numerical noise. Thus, one can observe ringing artifacts in the solution.

For $\sigma\in[0,\frac{1}{2})$ we have  $\re \psi_m(\sigma,\xi) \le 0$, i.e., we have no numerical damping. That is numerical noise
is not damped (for positive $\sigma$ there is even an amplification).
This effect usually causes convergence problems of the applied numerical methods, as e.g.\ 
Newton's method for nonlinear problems. In practice stable behavior was observed for the range $\sigma\in[-0.3,-0.1]$.
However, for oscillators we have the problem of numerical damping in this range.   

\section{Trigonometric splines for elimination of numerical damping \label{sec:trig}}

For the simulation of autonomous systems we need a method, with numerical damping at high frequencies, while low frequencies (in particular
the fundamental) frequency are not damped at all. This can be achieved if the numerical differentiation works exactly 
for the low frequencies (cf.~\cite{Gau61,Lampe2001,Bra2001}).
This can be achieved if splines are replaced by function spaces, which contain the low frequencies.
In order to preserve the useful properties of spline functions, trigonometric splines seem to be an interesting choice.
  
For an exhaustive description of trigonometric splines, even on non-uniform grids, we refer to \cite[Sect.~10.8]{Sch81}.
The space of trigonometric splines of order $m$ and mesh size $h$ is given as
$$
\widetilde{S}_{m,h}=\Big\{f\in C^{m-2}(\RR):~f\big|_{\big(h k,h(k+1)\big)}\in T_m\Big\},
$$
with the space of trigonometric polynomials
$$
T_m = \bigg\{\sum_{k=0}^m c_k\, e^{2\pi i (k-\frac{m-1}{2})t}:~c_k\in\CC\bigg\}.
$$
Obviously, for even order $T_{2\mu}$ is a space of anti-periodic functions, i.e.\  $f(t)=-f(t+1)$, and is not suited for our purpose.
In particular, constants and the fundamental frequency are not contained in $T_{2\mu}$ and $\widetilde{S}_{2\mu,h}$.
For odd order the space $T_{2\mu+1}$ contains the real valued trigonometric polynomials
$$
\sum_{k=0}^\mu a_k \cos(2\pi k t) + \sum_{k=1}^\mu a_k \sin(2\pi k t),
$$
which are therefore also contained in $\widetilde{S}_{2\mu+1,h}$.
A stable basis for $\widetilde{S}_{m,h}$ (if $h\, m<1$) is given by the translates $Q_{m,h}(t-hk)$ of the trigonometric B-spline $Q_{m,h}(t)$
defined by the recursion  
\begin{align*}
&Q_{1,h}(t) := \chi_{(0,h]}(t),\\[1ex]
&Q_{m,h}(t):=\\
&\quad\frac{\sin(\pi t)\,Q_{m-1,h}(t)+\sin\big(\pi (h\, m-t)\big)\,Q_{m-1,h}(t-h)}{\sin\big(\pi h (m-1)\big)}.
\end{align*}
One can see easily that $Q_{m,h}(t)$ is supported on $[0,m\,h]$ and from 
$\sin(t)=t+\mathcal{O}(t^3)$ we conclude
$Q_{m,h}(t)=N_m(ht)+\mathcal{O}(h^3)$, i.e., 
on fine grids the trigonometric splines behave similar to the
classical polynomial splines. 
The above formulation allows us also to take advantage of several spline algorithms
with only moderate extra computational effort for the computation of the sine function.

Analogously to  \ref{spline_diff} the numerical differentiation in a collocation method can be described by
$$
\widehat{y^\prime_k}=\widetilde\psi_{m,\frac{1}{n}}(\sigma,\tfrac{k}{n})\;\hat{y}_k,
$$
where 
$$ 
\widetilde\psi_{m,h}(x,\xi):=\frac{\tfrac{\partial}{\partial x}\widetilde\phi_{m,h}(x,\xi)}{\widetilde\phi_{m,h}(x,\xi)}
$$
and $\widetilde\phi_{m,h}(x,\xi)=\sum_k Q_{m,h}\big(h(x+\tfrac{m}{2}+k)\big)\,e^{2\pi i k \xi}$.

However, for $m=2\mu+1$ the interpolation of $e^{2\pi i k t}$, $k\in\ZZ$, $|k|<\mu$ is exact,
i.e., there are uniquely determined coefficients $c_{k,\ell}\in\CC$ such that
$$
e^{2\pi i k t}=s_k(t):=\sum_\ell c_{k,\ell}\, Q_{m,\frac{1}{n}}\big(\tfrac{t-\ell}{n}\big),\qquad|k|<\mu.
$$   
This implies in turn that the derivatives satisfy $s^\prime_k(t) = 2\pi i k\, s_k(t)$ or
$\widetilde\psi_{m,\frac{1}{n}}(\sigma,\tfrac{k}{n})=2\pi i k$, $|k|<\mu$.
That is, there is no damping of low frequencies (for $m=3,5,\ldots$). 

\begin{figure}[h]
	\center
	\includegraphics[width=\columnwidth]{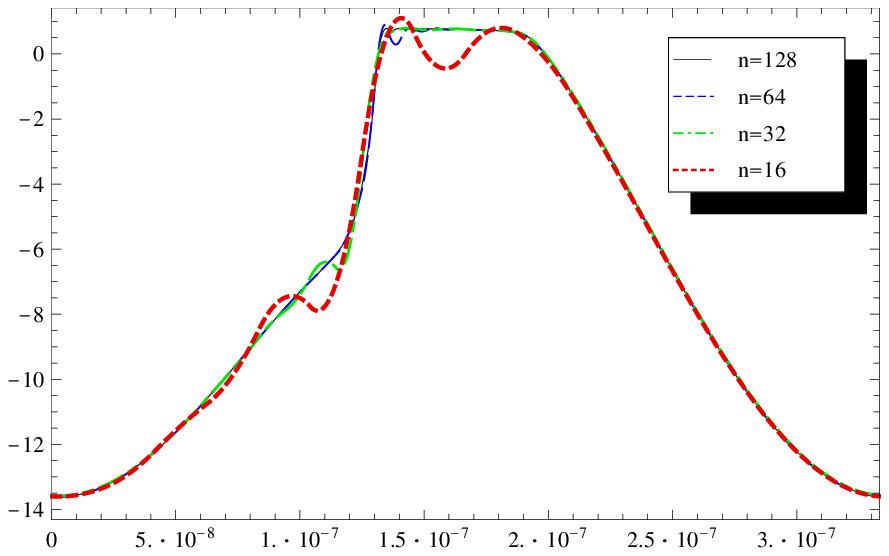}
	\fcaption{Trigonometric spline solutions for the 3MHz  Colpitz-Quartz oscillator, with $\sigma=-\frac{1}{4}$.}\label{f:fig5}
\end{figure}

Figure~\ref{f:fig5} shows the results of a simulation using trigonometric splines of order $m=3$,
with the same parameters as for  Fig.~\ref{f:fig1}. The numerical damping of the fundamental frequency
is eliminated, and the results behave in the range of the usual approximation error for the chosen grid.

\section{Conclusion}

We have shown that using trigonometric splines instead of the classical 
polynomial splines can eliminate the
negative effects of numerical damping for oscillator simulation, while we can still take advantage of many 
useful properties of B-splines. 

The methods can easily modified to nonuniform grids used for adaptive methods
described in \cite{BiBra12,BiDau10b,BiDau10a}. Here we may have grids 
with a locally high resolution in one
area and low resolution in other areas. 
Thus several 
negative effects studied in this article may occur simultaneously. 
In particular,
for adaptive grid refinement the differences in numerical damping between 
different grids can have serious
effects on the performance of the method, which do not occur if trigonometric 
splines are used.
Thus, the introduction of trigonometric splines is an important contribution to 
oscillator simulation. 

\Acknowledgements
This work was founded by the Austrian Science Fund (FWF): P22549-N18.


\begin{thebibliography}{10}

\bibitem{BiBra12}
K.~Bittner and H.-G. Brachtendorf.
\newblock Adaptive multirate wavelet method for circuit simulation.
\newblock In {\em 49th Design Automation Conference}, submitted.

\bibitem{BiDau10b}
K.~Bittner and E.~Dautbegovic.
\newblock Adaptiv wavelet-based method for simulation of electronic circuits.
\newblock In {\em Scientific Computing in Electrical Engineering 2010},
  Mathematics in Industry. Springer, Berlin Heidelberg, accepted.

\bibitem{BiDau10a}
K.~Bittner and E.~Dautbegovic.
\newblock Wavelets algorithm for circuit simulation.
\newblock In {\em Progress in Industrial Mathematics at ECMI 2010}, Mathematics
  in Industry. Springer, Berlin Heidelberg, submitted.

\bibitem{Bra2001}
H.~G. Brachtendorf.
\newblock {Theorie und Analyse von autonomen und quasiperiodisch angeregten
  elektrischen Netzwerken. Eine algorithmisch orientierte Betrachtung}.
\newblock Universit\"at Bremen, 2001.
\newblock {H}abilitationsschrift.

\bibitem{Bra2009}
H.~G. Brachtendorf, A.~Bunse-Gerstner, B.~Lang, and S.~Lampe.
\newblock Steady state of electronic circuits by cubic and exponential splines.
\newblock {\em Electrical Engineering}, 91:287--299, 2009.

\bibitem{Dau05}
E.~Dautbegovic.
\newblock {\em {Transient Simulation of Complex Electronic Circuits and Systems
  Operating at Ultra High Frequencies}}.
\newblock PhD thesis, Dublin City University, 2005.

\bibitem{DCB05}
E.~Dautbegovic, M.~Condon, and C.~Brennan.
\newblock An efficient nonlinear circuit simulation technique.
\newblock {\em IEEE Trans. Microwave Theory Tech.}, 53(2):548 --555, 2005.

\bibitem{Gau61}
W.~Gautschi.
\newblock {Numerical integration of ordinary differential equations based on
  trigonometric polynomials}.
\newblock {\em Numerische Mathematik}, 3:381--397, 1961.

\bibitem{JRS91}
K.~Jetter, S.~D. Riemenschneider, and N.~Sivakumar.
\newblock Schoenberg's exponential {Euler} spline.
\newblock {\em Proc.\ Royal Soc. Edinburgh}, 118A:21--35, 1991.

\bibitem{Lampe2001}
S.~Lampe, H.~G. Brachtendorf, E.~J.~W. {ter Maten}, S.~P. Onneweer, and
  R.~Laur.
\newblock Robust limit cycle calculations of oscillators.
\newblock pages 233--240, 2001.

\bibitem{Sch81}
L.~L. Schumaker.
\newblock {\em Spline Functions: Basic Theory}.
\newblock Wiley, New York, 1981.

\bibitem{SGN07}
N.~Soveiko, E.~Gad, and M.~Nakhla.
\newblock A wavelet-based approach for steady-state analysis of nonlinear
  circuits with widely separated time scales.
\newblock {\em IEEE Microwave Wireless Compon. Lett.}, 17(6):451--453, 2007.

\bibitem{SN03}
N.~Soveiko and M.~Nakhla.
\newblock Wavelet harmonic balance.
\newblock {\em Microwave and Wireless Components Letters, IEEE},
  13(6):232--234, 2003.

\end{thebibliography}


\end{document}